\documentclass{amsart}
\usepackage{lastpage}
\usepackage{yhmath,verbatim,mathabx}
\usepackage[all]{xy}

\DeclareFontFamily{U}{mathx}{\hyphenchar\font45}
\DeclareFontShape{U}{mathx}{m}{n}{
	<5> <6> <7> <8> <9> <10>
	<10.95> <12> <14.4> <17.28> <20.74> <24.88>
	mathx10
}{}
\DeclareSymbolFont{mathx}{U}{mathx}{m}{n}
\DeclareFontSubstitution{U}{mathx}{m}{n}
\DeclareMathAccent{\widecheck}{0}{mathx}{"71}
\DeclareMathAccent{\wideparen}{0}{mathx}{"75}

\newcommand{\bdis}{\begin{displaymath}}
\newcommand{\edis}{\end{displaymath}}
\newcommand{\be}{\begin{equation}}
\newcommand{\ee}{\end{equation}}
\newcommand{\mbb}{\mathbb}
\newcommand{\mcal}{\mathcal}

\newcommand{\vp}{\varphi}

\newcommand{\vth}{\vartheta}

\newcommand{\zf}{\zeta\left(\frac{1}{2}+it\right)}

\theoremstyle{definition}

\theoremstyle{remark}
\newtheorem{remark}[]{Remark}

\newtheorem*{mydef11}{{\bf Theorem 1}}

\newtheorem*{mydef12}{{\bf Theorem 2}}

\newtheorem*{mydef2}{{\bf Definition}}

\newtheorem*{mydef21}{{\bf Definition 1}}

\newtheorem*{mydef22}{{\bf Definition 2}}

\newtheorem*{mydef41}{{\bf Corollary 1}}

\newtheorem*{mydef42}{{\bf Corollary 2}}

\newtheorem*{mydef43}{{\bf Corollary 3}}

\newtheorem*{mydef51}{{\bf Lemma 1}}

\newtheorem*{mydef52}{{\bf Lemma 2}}

\newtheorem*{mydef53}{{\bf Lemma 3}}

\newtheorem*{mydef54}{{\bf Lemma 4}}

\newtheorem*{mydef8}{{\bf Property}}

\newtheorem*{mydef10}{{\bf Question}} 

\newtheorem*{mydef110}{{\bf Answer}}

\numberwithin{equation}{section}

\begin{document}

\title{Jacob's ladders, crossbreeding, secondary crossbreeding and synergetic phenomena generated by Riemann's zeta-function and some elementary functions on disconnected sets of the critical line}

\author{Jan Moser}

\address{Department of Mathematical Analysis and Numerical Mathematics, Comenius University, Mlynska Dolina M105, 842 48 Bratislava, SLOVAKIA}

\email{jan.mozer@fmph.uniba.sk}


\begin{abstract}
In this paper we obtain new complete hybrid formula for corresponding class of $\zeta$-factorization formulas. We demonstrate that this formula is the synergetic one. Namely, this one describes the cooperation between some class of elementary functions and Riemann's zeta-function on a class of disconnected sets on the critical line. Some comparizon between this pure mathematical phenomenon and the Belousov-Zhabotinski chemical reaction is added.  

\centering DEDICATED TO 110th ANNIVERSARY OF HARDY'S \\ PURE MATHEMATICS
\end{abstract}
\maketitle

\section{Introduction}

In this section we give: 
\begin{itemize} 
	\item[1.1] some notion about synergetics (cooperative) phenomena in systems that occurs in science, 
	\item[1.2] the simplest form of our main result, 
	\item[1.3] an interpretation of this result from the point of view of subsect. 1.1,  
	\item[1.4] final remark. 
\end{itemize}   

\subsection{} 

In sciences (physics, chemistry,biology, \dots) there are known various effects called as synergetic (cooperative) phenomena, (see \cite{1}). Within this context it is usually supposed that: 
\begin{itemize} 
\item[(a)] the system contain big number of objects (atoms, molecules or more complicated formations), 
\item[(b)] the cooperative phenomenon is such one in which the constituents of the system cannot be regarded as acting independently from each other, 
\item[(c)] cooperative phenomena result from interactions between constituents. 
\end{itemize} 

An wonderful example of such a phenomenon is given by the classical Belousov-Zhabotinski periodic chemical reaction (see \cite{1}, pp. 25-27). Namely 
\bdis 
\mbox{red} (excess \ of \ Ce^{3+})\longleftrightarrow \mbox{blue} (excess \ of \ Ce^{4+}) 
\edis  
oscillations of the chemical system on two disconnected parts of a test-tube (see \cite{1}, p. 27, Fig. 1.19). 

\subsection{} 

Now, we choose the following four elementary functions 
\be \label{1.1} 
\begin{split}
& (t-\pi L)\sin^2t, (t-\pi L)\cos^2t,\ (t-\pi L)^2\sin^2t,\ (t-\pi L)^2\cos^2t, \\ 
& t\in [\pi L,\pi L+U],\ U\in (0,\pi/2),\ L\in\mbb{N},\ L>L_0>0, 
\end{split}
\ee  
and high-transcendental function 
\be \label{1.2} 
\left|\zf\right|,\ t\in [\pi L_0,+\infty), 
\ee  
where $L_0$ is a big constant. 

The simplest form of our main result is the following one: for every sufficiently big $L$ and every $U\in (0,\pi/2)$ there are: 
\begin{itemize} 
\item[(a)] segment 
\bdis 
[\overset{1}{\wideparen{\pi L}},\overset{1}{\wideparen{\pi L+U}}]:\ [\pi L,\pi L+U]\prec [\overset{1}{\wideparen{\pi L}},\overset{1}{\wideparen{\pi L+U}}], 
\edis  
where the distance of the above mentioned segments $\rho(L)$ obeys the following 
\be \label{1.3} 
\rho(L)=\overset{1}{\wideparen{\pi L}}-(\pi L+U)\sim \pi (1-c)\pi(L)\sim \pi(1-c)\frac{L}{\ln L}\to\infty 
\ee  
as $L\to\infty$ (where $\pi(L)$ is the prime-counting function and $c$ is the Euler's constant), 
\item[(b)] the values 
\be \label{1.4} 
\begin{split} 
& x_1,\dots,x_4\in (\pi L,\pi L+U);\ x_1=x_1(U,L), \dots \\ 
& y_1,\dots,y_4,z\in 
(\overset{1}{\wideparen{\pi L}},\overset{1}{\wideparen{\pi L+U}}),\ 
y_1=y_1(U,L),\dots,\ z=z(U,L) 
\end{split} 
\ee  
such that the following formula 
\be \label{1.5} 
\begin{split} 
& \left\{ (x_1-\pi L)\sin^2 x_1
\left|\zeta\left(\frac 12+iy_1\right)\right|^2+
(x_2-\pi L)\cos^2 x_2
\left|\zeta\left(\frac 12+iy_2\right)\right|^2 \right\}^2 \\ 
& \sim 
\frac 34 \left|\zeta\left(\frac 12+iz\right)\right|^2\times \\ 
& \times 
\left\{ 
(x_3-\pi L)^2\sin^2 x_3\left|\zeta\left(\frac 12+iy_3\right)\right|^2+ 
(x_4-\pi L)^2\cos^2 x_4\left|\zeta\left(\frac 12+iy_4\right)\right|^2
\right\}  \\ 
& L\to\infty 
\end{split} 
\ee  
holds true on disconnected set 
\be \label{1.6} 
\Delta_1=[\pi L,\pi L+U]\bigcup [\overset{1}{\wideparen{\pi L}},\overset{1}{\wideparen{\pi L+U}}].
\ee 
\end{itemize}  

\begin{remark}
Let (see (\ref{1.4}), (\ref{1.5})) 
\bdis 
\begin{split} 
& A(x_1,x_2,y_1,y_2)= \\ 
& = \left\{ (x_1-\pi L)\sin^2 x_1
\left|\zeta\left(\frac 12+iy_1\right)\right|^2+
(x_2-\pi L)\cos^2 x_2
\left|\zeta\left(\frac 12+iy_2\right)\right|^2 \right\}^2, \\ 
& B(x_3,x_4,y_3,y_4)= \\ 
& = 
\left|\zeta\left(\frac 12+iz\right)\right|^2\times \\ 
& \times 
\left\{ 
(x_3-\pi L)^2\sin^2 x_3\left|\zeta\left(\frac 12+iy_3\right)\right|^2+ 
(x_4-\pi L)^2\cos^2 x_4\left|\zeta\left(\frac 12+iy_4\right)\right|^2
\right\} . 
\end{split}
\edis   
Then we can rewrite the result (\ref{1.5}) as the following law of 3/4: 
\be \label{1.7} 
\frac{A(x_1,x_2,y_1,y_2)}{B(x_3,x_4,y_3,y_4)}\sim \frac 34,\ L\to\infty . 
\ee 
\end{remark} 

\begin{remark}
Let us notice explicitly that the formula (\ref{1.5}), valid on disconnected set (\ref{1.6}), and obeying properties (\ref{1.3}) and (\ref{1.4}): 
\begin{itemize}
	\item[(a)] represents new type of result in the theory of the Riemann's zeta-function as well as in the theory of elementary functions 
	\item[(b)] it is impossible to achieve this one by the current methods in the theory of Riemann's zeta-function.  
\end{itemize} 
\end{remark} 

\subsection{} 

Let us notice the following: 
\begin{itemize}
	\item[(A)] We assign the corresponding set of values 
\be \label{1.8} 
\begin{split}  
& \{(t-\pi L)\sin^2t\},\ \{(t-\pi L)\cos^2t\},\ \{(t-\pi L)^2\sin^2t\}, \\ 
& \{(t-\pi L)^2\cos^2t\},\ t\in [\pi L,\pi L+U], \\ 
& \left\{\left|\zf\right|^2\right\},\ t\in [\overset{1}{\wideparen{\pi L}},\overset{1}{\wideparen{\pi L+U}}] . 
\end{split} 
\ee  
to every function of the sets (\ref{1.1}) and (\ref{1.2}). Every of these continuum sets represents in our case \emph{a more complicated formation} (comp. sect. 1.1, (a)). 
\item[(B)] Next, the set 
\be \label{1.9} 
 \{(t-\pi L)\sin^2t\}\cup \dots\cup \{(t-\pi L)^2\cos^2t\}\cup 
 \left\{\left|\zf\right|^2\right\},
\ee 
represents the main system on the disconnected set (\ref{1.6}) and the set (\ref{1.8}) are constituents of the system (\ref{1.9}) (comp. sect. 1.1, (b)). 
\item[(C)] Interactions between constituents (\ref{1.8}) are controlled by the Jacob's ladder (as it we shall see in the following). 
\item[(D)] Mentioned interactions generate the synergic (cooperative) phenomenon, i.e. deterministic formula (\ref{1.5}) (or (\ref{1.7})) on disconnected set (\ref{1.6}) together with the fundamental property (\ref{1.3}). 
\end{itemize} 

\subsection{} 

Finally, we give the following 

\begin{remark} 
The formulation of all results and proofs in this paper are based on new notions and methods in the theory of Riemann's zeta-function we have introduced in series of 45 papers concerning Jacob's ladders. These can be found in arXiv [math.CA] starting with the paper \cite{2}. 
\end{remark}  

In this paper we use especially the following notions: Jacob's ladder, algorithm for generating $\zeta$-factorization formulas, complete hybrid formula. We give the short survey of these notions in the following sects. 2 and 3 of this paper. 

Also new notions are introduced in this paper, namely: secondary crossbreeding and secondary complete hybrid formula. 

Finally, we notice explicitly that the results of this paper constitute the generic complement to the Riemann's functional equation on the critical line (comp. \cite{7}). 

\section{Jacob's ladders and the basic disconnected set (the first set of notions)} 

\subsection{} 

Let us remind that the Jacob's ladder 
\bdis 
\vp_1(t)=\frac 12\vp(t) 
\edis 
was introduced in \cite{2} (see also \cite{3}), where the function $\vp(t)$ is arbitrary solution of the nonlinear integral equation (also introduced in \cite{2}) 
\bdis 
\int_0^{\mu[x(T)]}Z^2(t)e^{-\frac{2}{x(T)}t}{\rm d}t=\int_0^TZ^2(t){\rm d}t, 
\edis  
where, of course, 
\bdis 
\begin{split} 
	& Z(t)=e^{i\vth(t)}\zf , \\ 
	& \vth(t)=-\frac t2\ln\pi+\mbox{Im}\ln\Gamma\left(\frac 14+i\frac t2\right), 
\end{split} 
\edis   
where each admissible function $\mu(y)$ generates a solution 
\bdis 
y=\vp(T;\mu)=\vp(T);\ \mu(y)>7y\ln y. 
\edis  
We call the function $\vp_1(T)$ the Jacob's ladder as an analogue of the Jacob's dream in Chumash, Bereishis, 28:12. 

\begin{remark} 
By making use of the Jacob's ladders we have shown (see \cite{2}) that the classical Hardy-Littlewood integral (1918) 
\bdis 
\int_0^T\left|\zf\right|^2{\rm d}t 
\edis  
has - in addition to previously known Hardy-Littlewood expression (and other similar) possessing an unbounded error term at $T\to\infty$ - the following infinite set of almost exact representations 
\bdis 
\begin{split} 
& \int_0^T\left|\zf\right|^2{\rm d}t=\vp_1(T)\ln\{\vp_1(T)\}+ \\ 
& + (c-\ln 2\pi)\vp_1(T)+c_0+\mcal{O}\left(\frac{\ln T}{T}\right),\ T\to\infty, 
\end{split} 
\edis  
where $c$ is the Euler's constant and $c_0$ is the constant from the Titschmarsh-Kober-Atkinson formula. 
\end{remark}  

\subsection{} 

Next, we have obtained (see \cite{2}, (6.2)) the following formula 
\be \label{2.1}  
T-\vp_1(T)\sim (1-c)\pi(T);\ \pi(T)\sim \frac{T}{\ln T},\ T\to\infty , 
\ee 
where $\pi(T)$ is the prime-counting function. 

\begin{remark}
Consequently, the Jacob's ladder $\vp_1(T)$ can be viewed by our formula (\ref{2.1}) as an asymptotic complementary function to the function 
\bdis 
(1-c)\pi(T) 
\edis  
in the following sense 
\bdis 
\vp_1(T)+(1-c)\pi(T)\sim T,\ T\to\infty. 
\edis 
\end{remark}  

Since Jacob's ladder is exactly increasing function, we have defined the reversely iterated sequence 
\bdis 
\{\overset{k}{T}\}_{k=1}^{k_0},\ k_0\in\mbb{N}:\ \vp_1(\overset{k}{T}):=\overset{k-1}{T},\ 
\overset{0}{T}=T>T_0>0;\ \overset{k}{T}=\vp_1^{-k}(T), 
\edis  
where $T_0$ is a sufficiently big, (we fix $k_0$). Next, we have proved (see \cite{4}, (2.5)-(2.7), (5.1)-(5.13)) the following properties: if 
\be \label{2.2} 
U=o\left(\frac{T}{\ln T}\right), \ T\to \infty 
\ee  
then 
\be \label{2.4} 
\begin{split}
&  \vp_1\{[\overset{k}{T},\overset{k}{\wideparen{T+U}}]\}=[\overset{k-1}{T},\overset{k-1}{\wideparen{T+U}}], \\
& [\overset{0}{T},\overset{0}{\wideparen{T+U}}]=[T,T+U], 
\end{split}
\ee 
\be \label{2.4} 
|[\overset{k}{T},\overset{k}{\wideparen{T+U}}]|=\overset{k}{\wideparen{T+U}}-\overset{k}{T}=o\left(\frac{T}{\ln T}\right), 
\ee 
\be \label{2.5} 
|[\overset{k-1}{\wideparen{T+U}},\overset{k}{T}]|=\overset{k}{T}-\overset{k-1}{\wideparen{T+U}}\sim(1-c)\pi(T), 
\ee  
\be \label{2.6} 
\begin{split} 
& [T,T+U]\prec [\overset{1}{T},\overset{1}{\wideparen{T+U}}]\prec \dots \prec [\overset{k}{T},\overset{k}{\wideparen{T+U}}], \\ 
& k=1,\dots,k_0. 
\end{split} 
\ee  
From (\ref{2.2})-(\ref{2.6})  we obtain the following basic property of Jacob's ladder. 

\begin{mydef8}
For every segment 
\bdis 
[T,T+U],\ U=o\left(\frac{T}{\ln T}\right), \ T\to \infty 
\edis  
there is the following class of disconnected sets 
\be \label{2.7} 
\Delta(T,U,k)=\bigcup_{r=0}^k [\overset{r}{T},\overset{r}{\wideparen{T+U}}],\ 1\leq k\leq k_0
\ee 
generated by the Jacob's ladder $\vp_1(T)$. 
\end{mydef8} 

\begin{remark} 
Disconnected set (\ref{2.7}) has the following properties: 
\begin{itemize}
	\item[(a)] lengths of its components are (see (\ref{2.4})) 
	\be \label{2.8} 
	|[\overset{r}{T},\overset{r}{\wideparen{T+U}}]|\sim o\left(\frac{T}{\ln T}\right),\ \ T\to \infty , 
	\ee  
	\item[(b)] lengths of its adjacent segments are (see (\ref{2.5})) 
	\be \label{2.9} 
	|[\overset{r-1}{\wideparen{T+U}},\overset{r}{T}]|\sim (1-c)\frac{T}{\ln T},\ T\to\infty . 
	\ee 
\end{itemize} 
\end{remark}  

\begin{remark}
The asymptotic behavior of the components of the disconnected set (\ref{2.7}) is as follows: if $T\to\infty$ then these components recede unboundedly each from other and all together are receding to infinity. Hence, at $T\to\infty$ the set of components of (\ref{2.7}) behaves as one dimensional Friedman-Hubble expanding universe. 
\end{remark} 

\section{$\zeta$-factorization formula, crossbreeding and complete hybrid formula (the second set of notions)} 

\subsection{} 

Let 
\be \label{3.1}  
f(t)\in \tilde{C}_0[T,T+U],\ U\in(0,U_0),\ U_0=o\left(\frac{T}{\ln T}\right),\ T\to\infty, 
\ee  
where 
\bdis 
f(t)\in \tilde{C}_0[T,T+U] \Leftrightarrow f(t)\in C[T,T+U] \wedge f(t)\geq 0,\notequiv 0, 
\edis 
of course, 
\be \label{3.2} 
F(U,T;f)=\frac 1U\int_T^{T+U}f(t){\rm d}t>0. 
\ee  
Next, by application of the operator $\hat{H}$ (that we have introduced in \cite{9}, (3.6)) on the element (\ref{3.1}) we obtain the following vector-valued function 
\be \label{3.3}  
\hat{H}f(t)=(\alpha_0^k,\alpha_1^k,\dots,\alpha_k^k,\beta_1^k,\dots,\beta_k^k),\ 1\leq k\leq k_0 
\ee  
(we fix arbitrary $k_0$), where 
\be \label{3.4}  
\begin{split}
&	\alpha_r^k=\alpha_r^k(U,T,k;f),\ r=0,1,\dots,k, \\ 
& 	\beta_r^k=\beta_r^k(U,T,k),\ r=1,\dots,k, \\ 
& \alpha_0^k\in (T,T+U), \\ 
& \alpha_r^k,\beta_r^k\in (\overset{r}{T},\overset{r}{\wideparen{T+U}}),\ r=1,\dots,k , 
\end{split}
\ee 
and simultaneously we obtain by our algorithm (see the short survey of this in \cite{9}, (3.1) - (3.11)), also the following asymptotic form of the $\zeta$-factorization formula 
\be \label{3.5} 
\begin{split} 
& \prod_{r=1}^k
\left|
\frac{\zeta\left(\frac 12+i\alpha_r^k\right)}{\zeta\left(\frac 12+i\beta_r^k\right)}
\right|^2= 
\left\{1+\mcal{O}\left(\frac{\ln\ln T}{\ln T}\right)\right\}\frac{F(U,T;f)}{f(\alpha_0^k)}\sim \\ 
& \sim \frac{F(U,T;f)}{f(\alpha_0^k)},\ T\to\infty , 
\end{split} 
\ee  
(see (\ref{3.2}),(\ref{3.4}), comp. \cite{8}, (2.1)-(2.7)). 

\begin{remark}
Since 
\be \label{3.6} 
\begin{split} 
& \alpha_r^k=\vp_1^{k-r}(d),\ d=d(U,T,k;f),\ r=0,1,\dots,k, \\ 
& \beta_r^k=\vp_1^{k-r}(e),\ e=e(U,T,k),\ r=1,\dots,k, 
\end{split}
\ee 
(see \cite{4}, (6.3), \cite{5}, (6.4), \cite{6}, (4.7), \cite{9}, (3.4), (3.5)) holds true, wee see that the functions $\alpha,\beta$ are generated by Jacob's ladder $\vp_1(t)$. Of course, we have: since 
\bdis 
d,e\in (\overset{r}{T},\overset{r}{\wideparen{T+U}})=(\vp_1^{-k}(T),\vp_1^{-k}(T+U)), 
\edis  
then 
\bdis 
\begin{split}
& \alpha_r^k=\vp_1^{k-r}(d)\in (\vp_1^{k-r}[\vp_1^{-k}(T)],\vp_1^{k-r}[\vp_1^{-k}(T+U)])= \\ 
& =  (\overset{r}{T},\overset{r}{\wideparen{T+U}}),\ \overset{r}{T}=\vp_1^{-r}(T), 
\end{split}
\edis  
and similar result for $\beta_r^k$. 
\end{remark} 

\subsection{} 

Now we begin with the following finite set 
\be \label{3.7} 
f_m(t)\in\tilde{C}_0[T,T+U], m=1,\dots,M,\ M\in\mbb{N} 
\ee  
of functions (condition for $U_0$ in (\ref{3.1}) is fulfilled). Similarly to (\ref{3.3})-(\ref{3.5}) we have in this case the following formulas 
\be \label{3.8} 
\hat{H}f_m=(\alpha_0^{m,k_m},\alpha_1^{m,k_m},\dots,\alpha_{k_m}^{m,k_m},\beta_1^{k_m},\beta_{k_m}^{k_m}),\ 1\leq k_m\leq k_0,\ k_m\in\mbb{N}, 
\ee  
where 
\be \label{3.9}  
\begin{split}
	&	\alpha_r^{m,k_m}=\alpha_r(U,T,k;f_m),\ r=0,1,\dots,k_m, \\ 
	& 	\beta_r^{k_m}=\beta_r^{k_m}(U,T,k),\ r=1,\dots,k_m, \\ 
	& \alpha_0^{m,k_m}\in (T,T+U), \\ 
	& \alpha_r^{m,k_m},\beta_r^{k_m}\in (\overset{r}{T},\overset{r}{\wideparen{T+U}}),\ r=1,\dots,k_m , 
\end{split}
\ee 
and 
\be \label{3.10}  
\begin{split} 
	& \prod_{r=1}^k
	\left|
	\frac{\zeta\left(\frac 12+i\alpha_r^{m,k_m}\right)}{\zeta\left(\frac 12+i\beta_r^{k_m}\right)}
	\right|^2\sim \frac{F_m(U,T;f_m)}{f_m(\alpha_0^{m,k_m})},\ m=1,\dots,M,\ T\to\infty .
\end{split} 
\ee   

\subsection{} 

Further, we will suppose (see \cite{10}, comp. \cite{11}) that after a finite number of stages of corresponding crossbreedings in the set (\ref{3.10}) (every member of this set participates on the crossbreeding) - that is after a finite number of eliminations of the external functions 
\be \label{3.11} 
F_m(U,T;f_m)>0 
\ee  
from the set (\ref{3.10}) - we obtain the following asymptotic complete hybrid formula (i.e. the result of complete elimination of the elements $F_m$, more exactly, of the corresponding external variables $U,T$): 
\be \label{3.12} 
\mcal{F}
\{
f_1(\alpha_0^{1,k_1})\prod_{r=1}^{k_1}(\dots),\dots,
f_M(\alpha_0^{M,k_M})\prod_{r=1}^{k_M}(\dots))\sim 1,\ T\to\infty, 
\}
\ee  
where, of course, 
\bdis 
\prod_{r=1}^{k_m}(\dots)=\prod_{r=1}^{k_m}\left|
\frac{\zeta\left(\frac 12+i\alpha_r^{m,k_m}\right)}{\zeta\left(\frac 12+i\beta_r^{k_m}\right)}
\right|^2 ,\ m=1,\dots,M. 
\edis  

\begin{remark}
In (\ref{3.1}) we may put, of course, 
\be \label{3.13} 
[T,T+U]\to [L,L+U],\ [\pi L,\pi L+U],\ \dots,\ L\in \mbb{N}. 
\ee 
\end{remark} 

\begin{remark}
Le us notice explicitly that we may also put 
\be \label{3.14} 
U\in (0,U_0)\to (0,U_0-\epsilon],\ [\epsilon,U_0)\to [\epsilon,U_0-\epsilon] 
\ee  
in (\ref{3.1}) if it is needful for existence of the inequalities 
\be \label{3.15} 
0\leq g(\epsilon)\leq F_m(U,T;f_m)\leq h(\epsilon),\ 
\ee  
for some intervals in (\ref{3.14}).  
\end{remark}  

For example, if  the formula (in the context (\ref{3.12})) 
\be \label{3.16} 
\mcal{F}_1\{\dots\}\sim \{ f_m(\alpha_0^{m,k_m})\prod_{r=1}^{k_m}\}^\Delta,\ \Delta\in\mbb{R},\ \Delta\not=0
\ee  
is the penultimate one in the process of crossbreeding then we obtain by (\ref{3.11}) or (\ref{3.15}) (comp. (\ref{3.12})) 
\be \label{3.17} 
\frac{\mcal{F}_1\{\dots\}}{\{ f_m(\alpha_0^{m,k_m})\prod_{r=1}^{k_m}\}^\Delta}\sim 
1,\ T\to\infty. 
\ee 

\begin{remark} 
Of course, 
\bdis 
\mbox{\{(3.16)\}} \ \Leftrightarrow \ \mbox{\{(3.17)\}}
\edis 
that is the formula (\ref{3.16}) is the asymptotic complete hybrid formula also. Consequently, we may break our process of crossbreeding on a penultimate formula (\ref{3.16}). 
\end{remark} 

\section{Lemmas} 

In this section we give four lemmas by making use of our algorithm for generating the $\zeta$-factorization formulae (see \cite{9}, (3.1) - (3.11), comp. sect. 3.1). 

\subsection{}  

Since 
\bdis 
\begin{split} 
& \frac 1U\int_{\pi L}^{\pi L+U}(t-\pi L)\sin^2t{\rm d}t=\\ 
& = \frac 14 U-\frac 14\sin 2U+\frac 14\frac{\sin^2U}{U}, 
\end{split} 
\edis  
then the have the following. 

\begin{mydef51} 
For the function 
\be \label{4.1} 
f_1(t)=f_1(t;L)=(t-\pi L)\sin^2t\in\tilde{C}_0[\pi L,\pi L+U],\ U\in (0,\pi/2),\ L\in\mbb{N}
\ee  
there are vector-valued functions 
\be \label{4.2} 
(\alpha_0^{1,k_1},\alpha_1^{1,k_1},\dots,\alpha_{k_1}^{1,k_1},\beta_1^{k_1},\dots,\beta_{k_1}^{k_1}),\ 1\leq k_1\leq k_0,\ k_1,k_0\in\mbb{N} 
\ee 
(here we fix arbitrary $k_0$) such that the following $\zeta$-factorization formula 
\be \label{4.3} 
\begin{split} 
& \prod_{r=1}^{k_1}\left|
\frac{\zeta\left(\frac 12+i\alpha_r^{1,k_1}\right)}{\zeta\left(\frac 12+i\beta_r^{k_1}\right)}
\right|^2\sim 
\frac{1}{(\alpha_0^{1,k_1}-\pi L)\sin^2\alpha_0^{1,k_1}}\times \\ 
& \times 
\left\{
\frac 14 U-\frac 14\sin 2U+\frac 14\frac{\sin^2U}{U}
\right\} , \ L\to\infty , 
\end{split}  
\ee 
holds true, where  
\be \label{4.4}  
\begin{split}
	&	\alpha_r^{1,k_1}=\alpha_r(U,\pi L,k_1;f_1),\ r=0,1,\dots,k_1, \\ 
	& 	\beta_r^{k_1}=\beta_r(U,\pi L,k_1),\ r=1,\dots,k_1, \\ 
	& \alpha_0^{1,k_1}\in (\pi L,\pi L+U), \\ 
	& \alpha_r^{1,k_1},\beta_r^{k_1}\in (\overset{r}{\wideparen{\pi L}},\overset{r}{\wideparen{\pi L+U}}),\ r=1,\dots,k_1 , 
\end{split}
\ee 
and 
\bdis 
(\overset{r}{\wideparen{\pi L}},\overset{r}{\wideparen{\pi L+U}})
\edis 
is the $r$th reverse iteration of the interval $(\pi L,\pi L+U)$ by means of the Jacob's ladder. 
\end{mydef51}

\subsection{} 

Since 
\bdis 
\begin{split} 
	& \frac 1U\int_{\pi L}^{\pi L+U}(t-\pi L)\cos^2t{\rm d}t=\\ 
	& = \frac 14 U+\frac 14\sin 2U-\frac 14\frac{\sin^2U}{U}, 
\end{split} 
\edis  
then the have the following. 

\begin{mydef52} 
	For the function 
	\be \label{4.5} 
	f_2(t)=f_2(t;L)=(t-\pi L)\cos^2t\in\tilde{C}_0[\pi L,\pi L+U],\ U\in (0,\pi/2),\ L\in\mbb{N}
	\ee  
	there are vector-valued functions 
	\be \label{4.6} 
	(\alpha_0^{2,k_2},\alpha_1^{2,k_2},\dots,\alpha_{k_2}^{2,k_2},\beta_1^{k_2},\dots,\beta_{k_2}^{k_2}),\ 1\leq k_2\leq k_0,\ k_2,k_0\in\mbb{N} 
	\ee 
	(here we fix arbitrary $k_0$) such that the following $\zeta$-factorization formula 
	\be \label{4.7} 
	\begin{split} 
		& \prod_{r=1}^{k_2}\left|
		\frac{\zeta\left(\frac 12+i\alpha_r^{2,k_2}\right)}{\zeta\left(\frac 12+i\beta_r^{k_2}\right)}
		\right|^2\sim 
		\frac{1}{(\alpha_0^{1,k_1}-\pi L)\cos^2\alpha_0^{2,k_2}}\times \\ 
		& \times 
		\left\{
		\frac 14 U+\frac 14\sin 2U-\frac 14\frac{\sin^2U}{U}
		\right\} , \ L\to\infty , 
	\end{split}  
	\ee 
	holds true, where  
	\be \label{4.8}  
	\begin{split}
		&	\alpha_r^{2,k_2}=\alpha_r(U,\pi L,k_2;f_2),\ r=0,1,\dots,k_2, \\ 
		& 	\beta_r^{k_2}=\beta_r(U,\pi L,k_2),\ r=1,\dots,k_2, \\ 
		& \alpha_0^{2,k_2}\in (\pi L,\pi L+U), \\ 
		& \alpha_r^{2,k_2},\beta_r^{k_2}\in (\overset{r}{\wideparen{\pi L}},\overset{r}{\wideparen{\pi L+U}}),\ r=1,\dots,k_2 .
	\end{split}
	\ee 
\end{mydef52} 

\subsection{}  

Since 
\bdis 
\begin{split} 
	& \frac 1U\int_{\pi L}^{\pi L+U}(t-\pi L)^2\sin^2t{\rm d}t=\\ 
	& = \frac 16 U^2-\frac 14\left(U^2-\frac 12\right)\frac{\sin 2U}{U}-\frac 14\cos 2U, 
\end{split} 
\edis  
then the have the following. 

\begin{mydef53} 
	For the function 
	\be \label{4.9} 
	f_3(t)=f_3(t;L)=(t-\pi L)^2\sin^2t\in\tilde{C}_0[\pi L,\pi L+U],\ U\in (0,\pi/2),\ L\in\mbb{N}
	\ee  
	there are vector-valued functions 
	\be \label{4.10} 
	(\alpha_0^{3,k_3},\alpha_1^{3,k_3},\dots,\alpha_{k_3}^{3,k_3},\beta_1^{k_3},\dots,\beta_{k_3}^{k_3}),\ 1\leq k_3\leq k_0,\ k_3,k_0\in\mbb{N} 
	\ee 
	(here we fix arbitrary $k_0$) such that the following $\zeta$-factorization formula 
	\be \label{4.11} 
	\begin{split} 
		& \prod_{r=1}^{k_3}\left|
		\frac{\zeta\left(\frac 12+i\alpha_r^{3,k_3}\right)}{\zeta\left(\frac 12+i\beta_r^{k_3}\right)}
		\right|^2\sim 
		\frac{1}{(\alpha_0^{3,k_3}-\pi L)^2\sin^2\alpha_0^{3,k_3}}\times \\ 
		& \times 
		\left\{
		 \frac 16 U^2-\frac 14\left(U^2-\frac 12\right)\frac{\sin 2U}{U}-\frac 14\cos 2U
		\right\} , \ L\to\infty , 
	\end{split}  
	\ee 
	holds true, where  
	\be \label{4.12}  
	\begin{split}
		&	\alpha_r^{3,k_3}=\alpha_r(U,\pi L,k_3;f_3),\ r=0,1,\dots,k_3, \\ 
		& 	\beta_r^{k_3}=\beta_r(U,\pi L,k_3),\ r=1,\dots,k_3, \\ 
		& \alpha_0^{3,k_3}\in (\pi L,\pi L+U), \\ 
		& \alpha_r^{3,k_3},\beta_r^{k_3}\in (\overset{r}{\wideparen{\pi L}},\overset{r}{\wideparen{\pi L+U}}),\ r=1,\dots,k_3 .
	\end{split}
	\ee 
\end{mydef53}

\subsection{}  

Since 
\bdis 
\begin{split} 
	& \frac 1U\int_{\pi L}^{\pi L+U}(t-\pi L)^2\cos^2t{\rm d}t=\\ 
	& = \frac 16 U^2+\frac 14\left(U^2-\frac 12\right)\frac{\sin 2U}{U}+\frac 14\cos 2U, 
\end{split} 
\edis  
then the have the following. 

\begin{mydef54} 
	For the function 
	\be \label{4.13} 
	f_4(t)=f_4(t;L)=(t-\pi L)^2\cos^2t\in\tilde{C}_0[\pi L,\pi L+U],\ U\in (0,\pi/2),\ L\in\mbb{N}
	\ee  
	there are vector-valued functions 
	\be \label{4.14} 
	(\alpha_0^{4,k_4},\alpha_1^{4,k_4},\dots,\alpha_{k_4}^{4,k_4},\beta_1^{k_4},\dots,\beta_{k_4}^{k_4}),\ 1\leq k_4\leq k_0,\ k_4,k_0\in\mbb{N} 
	\ee 
	(here we fix arbitrary $k_0$) such that the following $\zeta$-factorization formula 
	\be \label{4.15} 
	\begin{split} 
		& \prod_{r=1}^{k_4}\left|
		\frac{\zeta\left(\frac 12+i\alpha_r^{4,k_4}\right)}{\zeta\left(\frac 12+i\beta_r^{k_4}\right)}
		\right|^2\sim 
		\frac{1}{(\alpha_0^{4,k_4}-\pi L)^2\cos^2\alpha_0^{4,k_4}}\times \\ 
		& \times 
		\left\{
		\frac 16 U^2+\frac 14\left(U^2-\frac 12\right)\frac{\sin 2U}{U}+\frac 14\cos 2U
		\right\} , \ L\to\infty , 
	\end{split}  
	\ee 
	holds true, where  
	\be \label{4.16}  
	\begin{split}
		&	\alpha_r^{4,k_4}=\alpha_r(U,\pi L,k_4;f_4),\ r=0,1,\dots,k_4, \\ 
		& 	\beta_r^{k_4}=\beta_r(U,\pi L,k_4),\ r=1,\dots,k_4, \\ 
		& \alpha_0^{4,k_4}\in (\pi L,\pi L+U), \\ 
		& \alpha_r^{4,k_4},\beta_r^{k_4}\in (\overset{r}{\wideparen{\pi L}},\overset{r}{\wideparen{\pi L+U}}),\ r=1,\dots,k_4 .
	\end{split}
	\ee  
\end{mydef54}  

\section{Theorem: result of the interaction between a set of elementary functions and the function $|\zf|^2$} 

\subsection{} 

As a result of crossbreeding in the class of $\zeta$-factorization formulas (\ref{4.3}), (\ref{4.7}), (\ref{4.11}) and (\ref{4.15}) we will obtain the following. 

\begin{mydef11}
Let 
\be \label{5.1} 
f_l(t)\in\tilde{C}_0[\pi L,\pi L+U],\ U\in(0,\pi/2),\ l=1,2,3,4,\ L\in\mbb{N}, 
\ee  
where 
\be \label{5.2} 
\begin{split}
	& f_1(t)=(t-\pi L)\sin^2t, \\ 
	& f_2(t)= (t-\pi L)\cos^2t,\\
	& f_3(t)= (t-\pi L)^2\sin^2t,\\
	& f_4(t)= (t-\pi L)^2\cos^2t . 
\end{split}
\ee  
Then for the class of disconnected sets 
\be \label{5.3} 
\begin{split} 
& \Delta(U,\pi L,k_l)=\bigcup_{r=0}^{k_l}
[\overset{r}{\wideparen{\pi L}},\overset{r}{\wideparen{\pi L+U}}],\ 1\leq k_l\leq k_0, \\ 
& \Delta(U,\pi L,k_l)\subset \Delta(U,\pi L,k_0),\ l=1,2,3,4 
\end{split} 
\ee  
(here we fix arbitrary $k_0$) generated by the Jacob's ladder $\vp_1(t)$ there is the complete hybrid formula 
\be \label{5.4} 
\begin{split}
& \left\{
(\alpha_0^{1,k_1}-\pi L)\sin^2\alpha_0^{1,k_1}
\prod_{r=1}^{k_1}\left|
\frac{\zeta\left(\frac 12+i\alpha_r^{1,k_1}\right)}{\zeta\left(\frac 12+i\beta_r^{k_1}\right)}
\right|^2+ \right. \\ 
& \left. + (\alpha_0^{2,k_2}-\pi L)\cos^2\alpha_0^{2,k_2}
\prod_{r=1}^{k_2}\left|
\frac{\zeta\left(\frac 12+i\alpha_r^{2,k_2}\right)}{\zeta\left(\frac 12+i\beta_r^{k_2}\right)}
\right|^2
\right\}^2\sim \\ 
& \sim \frac 34 \times 
\left\{
(\alpha_0^{3,k_3}-\pi L)^2\sin^2\alpha_0^{3,k_3}
\prod_{r=1}^{k_3}\left|
\frac{\zeta\left(\frac 12+i\alpha_r^{3,k_3}\right)}{\zeta\left(\frac 12+i\beta_r^{k_3}\right)}
\right|^2+ \right. \\ 
& \left. + (\alpha_0^{4,k_4}-\pi L)^2\cos^2\alpha_0^{4,k_4}
\prod_{r=1}^{k_4}\left|
\frac{\zeta\left(\frac 12+i\alpha_r^{4,k_4}\right)}{\zeta\left(\frac 12+i\beta_r^{k_4}\right)}
\right|^2
\right\},\ L\to\infty 
\end{split}
\ee 
where 
\be \label{5.5}  
\begin{split} 
& 	\alpha_r^{l,k_l}=\alpha_r(U,\pi L,k_l;f_l),\ r=0,1,\dots,k_l, \\ 
& 	\beta_r^{k_l}=\beta_r(U,\pi L,k_l),\ r=1,\dots,k_l, \\ 
& \alpha_0^{l,k_l}\in (\pi L,\pi L+U),\ l=1,2,3,4, \\ 
& \alpha_r^{l,k_l},\beta_r^{k_l}\in (\overset{r}{\wideparen{\pi L}},\overset{r}{\wideparen{\pi L+U}}),\ r=1,\dots,k_l, 
\end{split} 
\ee 
next 
\be \label{5.6} 
\begin{split} 
& \rho\{
 [\overset{r-1}{\wideparen{\pi L}},\overset{r-1}{\wideparen{\pi L+U}}],
 [\overset{r}{\wideparen{\pi L}},\overset{r}{\wideparen{\pi L+U}}]
\}=\overset{r}{\wideparen{\pi L}}-\overset{r-1}{\wideparen{\pi L+U}}\sim \\ 
& \sim (1-c)\pi(\pi L)\sim \pi (1-c)\frac{L}{\ln L}\to\infty,\ L\to\infty , \\ 
& r=1,\dots,k_l,\ l=1,2,3,4, 
\end{split} 
\ee  
and, finally, (comp. (\ref{3.6})) 
\be \label{5.7} 
\begin{split}
& \alpha_r^{l,k_l}=\vp_1^{k_l-r}[d(U,\pi L,k_l;f_l)],\ r=0,1,\dots,k_l, \\ 
& \beta_r^{k_l}=\vp_1^{k_l-r}[e(U,\pi L,k_l)],\ r=1,\dots,k_l, \\ 
& d,e\in 
(\overset{k_l}{\wideparen{\pi L}},\overset{k_l}{\wideparen{\pi L+U}})=
(\vp_1^{-k_l}(\pi L),\vp_1^{-k_l}(\pi L+U)),\\ 
& l=1,2,3,4. 
\end{split}
\ee 
\end{mydef11}

\subsection{} 

Now, we use our notions and results from sects. 2 and 3 together with lemmas from sect. 4 to complete the proof of Theorem 1. The first stage of crossbreeding in the set of the four $\zeta$-factorization formulas (\ref{4.3}), (\ref{4.7}), (\ref{4.11}) and (\ref{4.15}) is expressed by the following formulas: first, (\ref{4.3}) and (\ref{4.7}) result in 
\be \label{5.8} 
\begin{split}
& (\alpha_0^{1,k_1}-\pi L)\sin^2\alpha_0^{1,k_1}
\prod_{r=1}^{k_1}\left|
\frac{\zeta\left(\frac 12+i\alpha_r^{1,k_1}\right)}{\zeta\left(\frac 12+i\beta_r^{k_1}\right)}
\right|^2+  \\ 
& . + (\alpha_0^{2,k_2}-\pi L)\cos^2\alpha_0^{2,k_2}
\prod_{r=1}^{k_2}\left|
\frac{\zeta\left(\frac 12+i\alpha_r^{2,k_2}\right)}{\zeta\left(\frac 12+i\beta_r^{k_2}\right)}
\right|^2\sim \frac 12U,\ L\to\infty , 
\end{split} 
\ee  
and, secondly, (\ref{4.11}) together with (\ref{4.15}) result in 
\be \label{5.9} 
\begin{split}
& (\alpha_0^{3,k_3}-\pi L)^2\sin^2\alpha_0^{3,k_3}
\prod_{r=1}^{k_3}\left|
\frac{\zeta\left(\frac 12+i\alpha_r^{3,k_3}\right)}{\zeta\left(\frac 12+i\beta_r^{k_3}\right)}
\right|^2+  \\ 
&  + (\alpha_0^{4,k_4}-\pi L)^2\cos^2\alpha_0^{4,k_4}
\prod_{r=1}^{k_4}\left|
\frac{\zeta\left(\frac 12+i\alpha_r^{4,k_4}\right)}{\zeta\left(\frac 12+i\beta_r^{k_4}\right)}
\right|^2\sim \frac 13U^2,\ L\to\infty. 
\end{split} 
\ee  
In the second stage of crossbreeding we eliminate the external $U$ from the formulae (\ref{5.8}) and (\ref{5.9}), where the external $U$ means that one contained on the right-hand sides of the mentioned formulas. As a result, we obtain the complete hybrid formula (\ref{5.4}). 

\subsection{} 

Now we give some notions about the interpretation of this result. 

We have the following continuum set of values of \emph{the active} functions 
\be \label{5.10} 
\begin{split}
& \{(t-\pi L)\sin^2t\},\ \{(t-\pi L)\cos^2t\}, \\ 
& \{(t-\pi L)^2\sin^2t\},\ \{(t-\pi L)^2\cos^2t\}, \\ 
& \pi L\leq t\leq \pi L+U, 
\end{split} 
\ee 
and 
\be \label{5.11} 
\begin{split}
& \left\{\left|\zeta\left(\frac 12+iy_r^l\right)\right|^2\right\},\ 
y_r^l\in [\overset{r}{\wideparen{\pi L}},\overset{r}{\wideparen{\pi L+U}}],\ r=1,\dots,k_l, \\ 
& l=1,2,3,4,\ 1\leq k_l\leq k_0,  
\end{split}
\ee  
on disconnected sets 
\be \label{5.12}  
\begin{split}
& \Delta(U,\pi L,k_l)=\bigcup_{r=0}^{k_l}
[\overset{r}{\wideparen{\pi L}},\overset{r}{\wideparen{\pi L+U}}],\ l=1,2,3,4, \\ 
& [\overset{0}{\wideparen{\pi L}},\overset{0}{\wideparen{\pi L+U}}]=[\pi L,\pi L+U], 
\end{split}
\ee 
where, of course, 
\be \label{5.13} 
\begin{split}
& \bigcup_{l=1}^4\Delta(U,\pi L,k_l)=\bigcup_{r=0}^{\bar{k}}
[\overset{r}{\wideparen{\pi L}},\overset{r}{\wideparen{\pi L+U}}], \\ 
& \bar{k}=\max\{ k_1,k_2,k_3,k_4\},\ L\to\infty. 
\end{split}
\ee 

\begin{remark}
The interactions of continuum sets (\ref{5.10}), (\ref{5.11}) (=constituents) that are controlled by the Jacob's ladder $\vp_1(t)$ generate the synergetic (cooperative) phenomenon, that is the complete hybrid formula (\ref{5.4}) for our system (\ref{5.10}) and (\ref{5.11}). Consequently, we have obtained a finite complex of $(k_0)^4$ formulas on the disconnected sets (\ref{5.12}) or (\ref{5.13}) whose neighboring components are separated by gigantic distances (see (\ref{5.6})). 
\end{remark}  

\begin{remark}
We notice explicitly that Jacob's ladder $\vp_1(t)$ makes a kind of double selection in our theory. namely, it selects: 
\begin{itemize}
	\item[(a)] disconnected set (\ref{5.13}) from every interval 
	\be \label{5.14} 
	[\pi L,+\infty),\ L\geq L_0>0
	\ee  
	for every admissible $U$ and $L_0$ sufficiently big, 
	\item[(b)] set of arguments (\ref{5.5}) (see (\ref{5.7})) in the disconnected (\ref{5.13}). 
\end{itemize} 
Consequently, properties (a) and (b) of the Jacob's ladder $\vp_1(t)$ create the base for the synergic formula (\ref{5.4}). 
\end{remark}

\begin{remark}
Next, we give some formal comparizon between the \emph{red and blue} in the Belousov-Zhabotinski chemical reaction and the disconnected set (\ref{5.13}) in our theory. Namely, we notice the following: 
\begin{itemize}
	\item[(a)] in the Belousov-Zhobotinski chemical reaction we have closed system and, consequently, duration of these disconnected sets is necessarily finite, 
	\item[(b)] in our theory, we have the synergetic formula (\ref{5.4}) based on every disconnected set (\ref{5.13}): 
	\be \label{5.15} 
	\forall\ L\in (\pi L_0,+\infty):\ L\to \Delta(U,\pi L,\bar{k}), 
	\ee  
	i.e. in every \emph{moment} there is infinite class of disconnected sets $\Delta(U,\pi L,\bar{k})$ - since our system is not closed (see (\ref{5.10}), (\ref{5.11}) and (\ref{5.13})). 
\end{itemize}
\end{remark} 

\subsection{} 

Now we give the following 
\begin{mydef10}
What is a ghost that animates the class of the following sets 
\be\label{5.16}  
\{ f_1(t)\}, \{ f_2(t)\}, \{f_3(t)\}, \{f_4(t)\},\ \{|\zf|^2\} ,\ t\geq \pi L, 
\ee  
where, in this case, \emph{to animate} means to excite an interaction between the constituents (\ref{5.16})? 
\end{mydef10} 

In this direction we give the following 

\begin{remark}
These are the levels of the animation. 
\begin{itemize}
\item[(a)] In the beginning, we have non-reactive sets (\ref{5.16}) - a mixture of inert gases. 
\item[(b)] Next, we have the origin of $\vp_1(t)$: 
\bdis 
|\zf|^2\to \vp_1(t), 
\edis 
that is the function $|\zeta|^2$ generates new function, namely, the Jacob's ladder (see sect. 2.1). 
\item[(c)] Further, Jacob's ladder generates: 
\begin{itemize}
	\item[(1)] the set of iterations 
	\be \label{5.17} 
	\{\vp_1^r(t)\}_{r=0}^{k_l},\ l=1,2,3,4, 
	\ee  
	\item[(2)] the class of disconnected sets (comp. (\ref{2.7})) 
	\bdis 
	\Delta(U,\pi L,k_l),\ \forall\- L\geq L_0. 
	\edis 
\end{itemize} 
\item[(d)] Finally, we have the following compositions 
\be \label{5.18} 
\{f_l(\alpha_0^{l,k_l})\},\ 
\left\{ \left|\zeta\left(\frac 12+i\alpha_r^{l,k_l}\right)\right|^2 \right\},\ 
\left\{ \left|\zeta\left(\frac 12+i\beta_r^{k_l}\right)\right|^2 \right\},\ l=1,2,3,4 
\ee  
of the sets (\ref{5.16}) with the Jacob's ladder $\vp_1(t)$ and its iterations (\ref{5.17}) (see (\ref{5.7}) and Remark 8). 
\item[(e)] Consequently, synergetic formula (\ref{5.4}) is created from the set of values  (\ref{5.18}). 
\end{itemize} 
\end{remark} 

Now we have the following. 

\begin{mydef110}
It is just the Jacob's ladder that is the ghost animating the class of sets (\ref{5.16}). 
\end{mydef110} 

\section{Some consequences from Theorem 1} 

\subsection{} 

Let us remind (see \cite{10}, comp. \cite{11}) the following notion. 

\begin{mydef2}
The subset 
\bdis 
\{ f_1(t),\dots,f_M(t)\},\ t\in [T,T+U] 
\edis  
of the set of (admissible) functions for which there is the complete hybrid formula is called as the family of $\zeta$-kindred functions, or shortly as the $\zeta$-family of elements. 
\end{mydef2}  

As a consequence we have by (\ref{5.2}), (\ref{5.4}) the following 

\begin{mydef41}
The subset 
\be \label{6.1} 
\begin{split} 
& \{(t-\pi L)\sin^2t, (t-\pi L)\cos^2t,\ (t-\pi L)^2\sin^2t,\ (t-\pi L)^2\cos^2t\}, \\ 
& t\in [\pi L,\pi L+U],\ U\in (0,\pi/2),\ L\in\mbb{N}, 
\end{split}
\ee 
is the $\zeta$-family of elements in the class $\tilde{C}_0[\pi L,\pi L+U]$. 
\end{mydef41} 

\subsection{} 

Next, in the case 
\be \label{6.2} 
k_1=k_2=k_3=k_4=k,\ 1\leq k\leq k_0 
\ee  
we obtain from Theorem 1 the following 

\begin{mydef42}
\be \label{6.3}
\begin{split}
& \left\{
(\alpha_0^{1,k}-\pi L)\sin^2\alpha_0^{1,k}
\prod_{r=1}^{k}\left|
\frac{\zeta\left(\frac 12+i\alpha_r^{1,k}\right)}{\zeta\left(\frac 12+i\beta_r^{k}\right)}
\right|^2+ \right. \\ 
& \left. + (\alpha_0^{2,k}-\pi L)\cos^2\alpha_0^{2,k}
\prod_{r=1}^{k}\left|
\frac{\zeta\left(\frac 12+i\alpha_r^{2,k}\right)}{\zeta\left(\frac 12+i\beta_r^{k}\right)}
\right|^2
\right\}^2\sim \\ 
& \sim \frac 34 \times 
\left\{
(\alpha_0^{3,k}-\pi L)^2\sin^2\alpha_0^{3,k}
\prod_{r=1}^{k}\left|
\frac{\zeta\left(\frac 12+i\alpha_r^{3,k}\right)}{\zeta\left(\frac 12+i\beta_r^{k}\right)}
\right|^2+ \right. \\ 
& \left. + (\alpha_0^{4,k}-\pi L)^2\cos^2\alpha_0^{4,k}
\prod_{r=1}^{k}\left|
\frac{\zeta\left(\frac 12+i\alpha_r^{4,k}\right)}{\zeta\left(\frac 12+i\beta_r^{k}\right)}
\right|^2
\right\}
\end{split}
\ee 
where 
\bdis 
\alpha_0^{l,k}\in (\pi L,\pi L+U),\ \alpha_r^{l,k},\beta_r^k\in 
(\overset{r}{\wideparen{\pi L}},\overset{r}{\wideparen{\pi L+U}}), 
\edis  
and 
\bdis 
r=1,\dots,k,\ l=1,2,3,4. 
\edis 
\end{mydef42} 

\begin{remark}
For example, in the case 
\bdis 
k=k_0=10^6 
\edis   
the basic disconnected set 
\bdis 
\Delta(U,\pi L,10^6) 
\edis  
has $10^6+1$ components and: 
\begin{itemize}
	\item[(a)] the component $[\pi L,\pi L+U]$ contains four arguments, 
	\item[(b)] every component 
	\bdis 
	[\overset{r}{\wideparen{\pi L}},\overset{r}{\wideparen{\pi L+U}}],\ r=1,\dots,10^6 
	\edis 
	contains five arguments. 
\end{itemize}
\end{remark} 

\begin{remark}
For comparison we notice that the main synergetic formula (\ref{5.4}) in the case $k_0=10^6$ contains $10^{24}$ morphologically different formulas. 
\end{remark} 

\begin{remark}
The case $k=k_0=1$ we have presented as the simplest result in the form (\ref{1.5}). 
\end{remark} 

\section{Secondary crossbreeding - Definition} 

\subsection{} 

We have defined the following notions in the sect. 3 of this paper: 
\begin{itemize}
	\item[(a)] our algorithm for generating the $\zeta$-factorization formula for a given admissible real function, 
	\item[(b)] operation of crossbreeding on a given class of $\zeta$-factorization formulas together with the notion of the complete hybrid formula (asymptotic, say) as the final product of this crossbreeding. 
\end{itemize} 
Now, in this section of our paper, we introduce another notion, namely the notion of secondary crossbreeding for every two admissible complete hybrid formulae. 

\subsection{} 

We begin with the set of two complete hybrid formulae (comp. (\ref{3.7}) - (\ref{3.12})): 
\be \label{7.1} 
\begin{split}
&\mcal{F}_1
\left\{
f_1(\alpha_0^{1,k_1})
\prod_{r=1}^{k_1}\left|
\frac{\zeta\left(\frac 12+i\alpha_r^{1,k_1}\right)}{\zeta\left(\frac 12+i\beta_r^{k_1}\right)}
\right|^2,\dots, \right. \\ 
& \left. 
f_M(\alpha_0^{M,k_M})
\prod_{r=1}^{k_M}\left|
\frac{\zeta\left(\frac 12+i\alpha_r^{M,k_M}\right)}{\zeta\left(\frac 12+i\beta_r^{k_M}\right)}
\right|^2
\right\} \sim 1\\ 
&  1\leq k_1,\dots,k_M\leq k_0,\ L\to\infty, 
\end{split}
\ee  
and 
\be \label{7.2} 
\begin{split}
	&\mcal{G}_1
	\left\{
	f_{M+1}(\alpha_0^{M+1,k_{M+1}})
	\prod_{r=1}^{k_{M+1}}\left|
	\frac{\zeta\left(\frac 12+i\alpha_r^{{M+1},k_{M+1}}\right)}{\zeta\left(\frac 12+i\beta_r^{k_{M+1}}\right)}
	\right|^2,\dots, \right. \\ 
	& \left. 
	f_{M+P}(\alpha_0^{{M+P},k_{M+P}})
	\prod_{r=1}^{k_{M+P}}\left|
	\frac{\zeta\left(\frac 12+i\alpha_r^{M+P,k_{M+P}}\right)}{\zeta\left(\frac 12+i\beta_r^{k_{M+P}}\right)}
	\right|^2
	\right\}\sim 1 , \\ 
&	1\leq k_{M+1},\dots,k_{M+P}\leq k_0,\ L\to\infty.  
\end{split}
\ee 
Next, we put, for example, in (\ref{7.1}), 
\bdis 
k_1=k_2=\dots=k_M=k,\ 1\leq k\leq k_0, 
\edis  
and we will assume that for the resulting formula there is an expression (an asymptotic solution) 
\be \label{7.3} 
\begin{split}
& \prod_{r=1}^k\left|\zeta\left(\frac 12+i\beta_r^k\right)\right|^2\sim \\ 
& \sim 
\mcal{F}_2
\left\{
f_1(\alpha_0^{1,k})\prod_{r=1}^k\left|\zeta\left(\frac 12+i\alpha_r^{1,k}\right)\right|^2, \dots , 
f_M(\alpha_0^{M,k})\prod_{r=1}^k\left|\zeta\left(\frac 12+i\alpha_r^{M,k}\right)\right|^2
\right\} , \\ 
& 1\leq k\leq k_0,\ L\to\infty. 
\end{split} 
\ee  

\begin{remark}
Of course, the $\beta$\-'s are independent on the choice of 
\bdis  
f_l(t),\ l=1,\dots,M+P, 
\edis  
and, more generally, on the choice of every admissible $f(t)$. 
\end{remark} 

If we put into (\ref{7.3}) 
\bdis 
k=k_{M+1}=\dots=k_{M+P} 
\edis  
(comp. inequalities in (\ref{7.1}) - (\ref{7.3})), then we obtain following formulae 
\be \label{7.4} 
\begin{split}
& \prod_{r=1}^{k_{M+l}}\left|\zeta\left(\frac 12+i\beta_r^{k_M+l}\right)\right|^2= \\ 
& = \frac{1}{1+o_l(1)}\times \\ 
&  \times
\mcal{F}_2
\left\{
f_1(\alpha_0^{1,k_M+l})\prod_{r=1}^{k_M+l}\left|\zeta\left(\frac 12+i\alpha_r^{1,k_{M+l}}\right)\right|^2, \dots , \right. \\ 
& \left.  
f_M(\alpha_0^{M,k_{M+l}})\prod_{r=1}^{k_M+l}\left|\zeta\left(\frac 12+i\alpha_r^{M,k_{M+l}}\right)\right|^2
\right\}= \\ 
& = \frac{1}{1+o_l(1)}\mcal{F}_2^{k_{M+l}}\{\dots\} , \\ 
& l=1,\dots,P,\ L\to\infty . 
\end{split}
\ee  
Further, we put (\ref{7.4}) into (\ref{7.2}) that gives 
\be \label{7.5} 
\begin{split}
&\mcal{G}_1
\left\{
\{1+o_1(1)\}\frac{f_{M+1}(\alpha_0^{M+1,k_{M+1}})}{\mcal{F}_2^{k_{M+1}}\{\dots\}}
\prod_{r=1}^{k_{M+1}}\left|\zeta\left(\frac 12+i\alpha_r^{M+1,k_{M+l}}\right)\right|^2,\dots , 
\right. \\ 
& \left. 
\{1+o_P(1)\}\frac{f_{M+P}(\alpha_0^{M+P,k_{M+P}})}{\mcal{F}_2^{k_{M+P}}\{\dots\}}
\prod_{r=1}^{k_{M+P}}\left|\zeta\left(\frac 12+i\alpha_r^{M+P,k_{M+P}}\right)\right|^2
\right\} = \\ 
& = 1+o(1). 
\end{split}
\ee

Finally, we will assume that we can obtain the following one 
\be \label{7.6} 
\begin{split}
& \{1+o(1)\}\times \\ 
& \times 
\mcal{G}_1 
\left\{
\frac{f_{M+1}(\alpha_0^{M+1,k_{M+1}})}{\mcal{F}_2^{k_M+1}\{\dots\}}
\prod_{r=1}^{k_{M+1}}\left|
\frac{\zeta\left(\frac 12+i\alpha_r^{{M+1},k_{M+1}}\right)}{\zeta\left(\frac 12+i\beta_r^{k_{M+1}}\right)}
\right|^2,\dots, \right. \\ 
& \left. 
\frac{f_{M+P}(\alpha_0^{{M+P},k_{M+P}})}{\mcal{F}_2^{k_M{M+P}}\{\dots\}}
\prod_{r=1}^{k_{M+P}}\left|
\frac{\zeta\left(\frac 12+i\alpha_r^{M+P,k_{M+P}}\right)}{\zeta\left(\frac 12+i\beta_r^{k_{M+P}}\right)}
\right|^2
\right\}=1+o(1), \\ 
& \Longrightarrow  \\ 
& \mcal{G}_1\{\dots\}\sim 1,\ L\to\infty  
\end{split}
\ee 
for this formula. 

\subsection{} 

Now we give the following definition. 

\begin{mydef21}
In this case we will call: 
\begin{itemize}
	\item[(a)] the procedure of elimination (\ref{7.3}) - (\ref{7.6}) as the secondary crossbreeding between two asymptotic complete hybrid formulae (\ref{5.4}) and (\ref{7.1}), 
	\item[(b)] the final result (\ref{7.6}) of this procedure as the asymptotic secondary complete hybrid formula. 
\end{itemize}
\end{mydef21} 

\begin{mydef22}
We will call the set of corresponding active functions (in our process) 
\bdis 
\{f_1(t),\dots,f_M(t)\}\cup \{ f_{M+1}(t),\dots,f_{M+P}(t)\}= \{ f_1(t),\dots ,f_{M+P}(t)\}
\edis 
as the $\zeta$-family of the second order. 
\end{mydef22} 

\begin{remark}
Of course, the secondary complete hybrid formula (\ref{7.6}) represents also the synergetic phenomenon generated by the interactions between sets of values 
\bdis 
\{f_1(t)\},\ \{f_2(t)\},\ \dots,\ \{f_{M+P}(t)\},\ t\in [\pi L_0, +\infty) 
\edis  
of elementary (say) functions and the set of values 
\bdis 
\left\{\zf\right\},\ t\in [\pi L_0, +\infty) 
\edis 
of highly transcendental function on the class of basic disconnected sets (\ref{2.7}). 
\end{remark} 

\section{Secondary crossbreeding - Example} 

\subsection{} 

Let us remind that we have obtained in our paper \cite{10} (comp. \cite{8}) the following complete hybrid formula 
\be \label{8.1} 
\begin{split}
& \sin^2\alpha_0^{5,k_5}
\prod_{r=1}^{k_5}
\left|\frac{\zeta\left(\frac 12+i\alpha_r^{{5},k_{5}}\right)}{\zeta\left(\frac 12+i\beta_r^{k_{5}}\right)} 
\right|^2 + \\ 
& + 
\cos^2\alpha_0^{6,k_6}
\prod_{r=1}^{k_6}
\left|\frac{\zeta\left(\frac 12+i\alpha_r^{{6},k_{6}}\right)}{\zeta\left(\frac 12+i\beta_r^{k_{6}}\right)} 
\right|^2\sim 1,\ 1\leq k_5,k_6\leq k_0,\ L\to\infty, 
\end{split} 
\ee  
where the function 
\be \label{8.2} 
f_5(t)=\sin^2(t)\in\tilde{C}_0[\pi L,\pi L+U],\ U\in (0,\pi/2) 
\ee  
generates the vector-valued functions 
\be \label{8.3}  
\begin{split}
    & (\alpha_0^{5,k_5},\alpha_1^{5,k_5},\dots,\alpha_{k_5}^{5,k_5},\beta_1^{k_5},\dots,\beta_{k_5}^{k_5}),  \\ 
	&	\alpha_r^{5,k_5}=\alpha_r(U,L,k_5;f_5),\ r=0,1,\dots,k_5, \\ 
	& 	\beta_r^{k_5}=\beta_r(U,L,k_5),\ r=1,\dots,k_5, \\ 
	& \alpha_0^{5,k_5}\in (\pi L,\pi L+U), \\ 
	& \alpha_r^{5,k_5},\beta_r^{k_5}\in (\overset{r}{\wideparen{\pi L}},\overset{r}{\wideparen{\pi L+U}}),\ r=1,\dots,k_5 , \ 1\leq k_5\leq k_0,  
\end{split}
\ee  
and the function 
\be \label{8.4} 
f_6(t)=\cos^2(t)\in\tilde{C}_0[\pi L,\pi L+U],\ U\in (0,\pi/2) 
\ee  
generates its vector-valued functions 
\be \label{8.5}  
\begin{split}
	& (\alpha_0^{6,k_6},\alpha_1^{6,k_6},\dots,\alpha_{k_6}^{6,k_6},\beta_1^{k_6},\dots,\beta_{k_6}^{k_6}),  \\ 
	&	\alpha_r^{6,k_6}=\alpha_r(U,L,k_6;f_6),\ r=0,1,\dots,k_6, \\ 
	& 	\beta_r^{k_6}=\beta_r(U,L,k_6),\ r=1,\dots,k_6, \\ 
	& \alpha_0^{6,k_6}\in (\pi L,\pi L+U), \\ 
	& \alpha_r^{6,k_6},\beta_r^{k_6}\in (\overset{r}{\wideparen{\pi L}},\overset{r}{\wideparen{\pi L+U}}),\ r=1,\dots,k_6 , \ 1\leq k_6\leq k_0. 
\end{split}
\ee   

\subsection{} 

Here we give the secondary complete hybrid formula as a product of secondary crossbreeding between two complete hybrid formulae (\ref{5.4}) and (\ref{8.1}). 

If we put into (\ref{8.1}) 
\bdis 
k_5=k_6=k,\ 1\leq k\leq k_0 
\edis  
and after this consecutively 
\bdis 
k=k_1,\dots,k_4, 
\edis  
then we obtain the following set of formulas 
\be \label{8.6} 
\begin{split}
& \prod_{r=1}^{k_l}
\left|\zeta\left(\frac 12+i\beta_r^{k_{l}}\right)\right|^2= \frac{1}{1+o_l(1)}\times \\ 
& \times 
\left\{
\sin^2\alpha_0^{5,k_l}
\prod_{r=1}^{k_l}
\left|\zeta\left(\frac 12+i\alpha_r^{5,k_{l}}\right)\right|^2+
\cos^2\alpha_0^{6,k_l}
\prod_{r=1}^{k_l}
\left|\zeta\left(\frac 12+i\alpha_r^{6,k_{l}}\right)\right|^2
\right\} ,\\ 
& l=1,2,3,4. 
\end{split}
\ee 
Next, after the substitution of (\ref{8.6}), $l=1,2$, into the left-hand side of (\ref{5.4}) we obtain (by some small algebra, comp. (\ref{7.6})) 
\bdis 
\begin{split}
& \{ A_3\}^2:\ A_3=\{1+o_1(1)\}A_1+\{1+o_2(1)\}A_2= \\ 
& = (A_1+A_2)\left(1+o_1(1)\frac{A_1}{A_1+A_2}+o_2(1)\frac{A_1}{A_1+A_2}\right)= \\ 
& = (1+o_3(1))(A_1+A_2);\ A_1,A_2>0 , 
\end{split}
\edis  
and a similar result for the right-hand side of (\ref{5.4}). Consequently, the complete elimination of the corresponding products in (\ref{5.4}) by means of (\ref{8.6}) gives the following short version (comp. Theorem 1) of 

\begin{mydef12}
\be \label{8.7} 
\begin{split}
& \left\{
\frac
{(\alpha_0^{1,k_1}-\pi L)\sin^2\alpha_0^{1,k_1}\prod_{r=1}^{k_1}\left|\zeta\left(\frac 12+i\alpha_r^{1,k_{1}}\right)\right|^2}
{\sin^2\alpha_0^{5,k_1}
	\prod_{r=1}^{k_1}
	\left|\zeta\left(\frac 12+i\alpha_r^{5,k_{1}}\right)\right|^2+
	\cos^2\alpha_0^{6,k_1}
	\prod_{r=1}^{k_1}
	\left|\zeta\left(\frac 12+i\alpha_r^{6,k_{1}}\right)\right|^2}+ \right. \\ 
& \left. + 
\frac
{(\alpha_0^{2,k_2}-\pi L)\cos^2\alpha_0^{2,k_2}\prod_{r=1}^{k_2}\left|\zeta\left(\frac 12+i\alpha_r^{2,k_{2}}\right)\right|^2}
{\sin^2\alpha_0^{5,k_2}
	\prod_{r=1}^{k_2}
	\left|\zeta\left(\frac 12+i\alpha_r^{5,k_{2}}\right)\right|^2+
	\cos^2\alpha_0^{6,k_2}
	\prod_{r=1}^{k_2}
	\left|\zeta\left(\frac 12+i\alpha_r^{6,k_{2}}\right)\right|^2} 
\right\}^2 \\ 
& \sim \\ 
& \frac 34 \times \\ 
& \left\{
\frac
{(\alpha_0^{3,k_3}-\pi L)^2\sin^2\alpha_0^{3,k_3}\prod_{r=1}^{k_3}\left|\zeta\left(\frac 12+i\alpha_r^{3,k_{3}}\right)\right|^2}
{\sin^2\alpha_0^{5,k_3}
	\prod_{r=1}^{k_3}
	\left|\zeta\left(\frac 12+i\alpha_r^{5,k_{3}}\right)\right|^2+
	\cos^2\alpha_0^{6,k_3}
	\prod_{r=1}^{k_3}
	\left|\zeta\left(\frac 12+i\alpha_r^{6,k_{3}}\right)\right|^2}+ \right. \\ 
& \left. + 
\frac
{(\alpha_0^{4,k_4}-\pi L)^2\cos^2\alpha_0^{4,k_4}\prod_{r=1}^{k_4}\left|\zeta\left(\frac 12+i\alpha_r^{4,k_{4}}\right)\right|^2}
{\sin^2\alpha_0^{5,k_4}
	\prod_{r=1}^{k_4}
	\left|\zeta\left(\frac 12+i\alpha_r^{5,k_{4}}\right)\right|^2+
	\cos^2\alpha_0^{6,k_4}
	\prod_{r=1}^{k_4}
	\left|\zeta\left(\frac 12+i\alpha_r^{6,k_{4}}\right)\right|^2} 
\right\}, \\ 
& 1\leq k_1,\dots,k_4\leq k_0,\ L\to\infty. 
\end{split}
\ee 
\end{mydef12}

\begin{remark}
Of course, the secondary complete hybrid formula (\ref{8.7}) represents the synergetic phenomenon also (comp. Remark 20). 
\end{remark} 

We give also the following consequence from our Theorem 2 (comp. Definition 2) 
\begin{mydef43}
The set of functions 
\bdis 
\begin{split}
& \{\sin^2t,\cos^2t,(t-\pi L)\sin^2t,(t-\pi L)\cos^2t,(t-\pi L)^2\sin^2t,(t-\pi L)^2\cos^2t\}, \\ 
& t\in [\pi L,\pi L+U],\ U\in (0,\pi/2),\  L\to\infty 
\end{split}
\edis 
is the $\zeta$-family of the second order. 
\end{mydef43} 

\section{On complicated internal structure of formulas that we have presented in this paper} 

Every formula as $\zeta$-factorization formula, complete hybrid formula and secondary complete hybrid formula contains product of type 
\bdis 
\prod_{r=1}^{k}
\left|\zeta\left(\frac 12+i\alpha_r\right)\right|^2. 
\edis  
Consequently, it is sufficient  to demonstrate the complicated internal structure of 
\bdis 
|Z(t)|=\left|\zf\right| 
\edis  
(comp. sect. 2.1). For this purpose we use the spectral form of the Riemann-Siegel formula (see \cite{6}, (3.1) - (3.8)) 
\bdis 
\begin{split}
& Z(t)=\sum_{n\leq \tau(x)}\frac{2}{\sqrt{n}}\cos\{t\omega_n(x)+\psi(x)\}+R(x),\\ 
& \tau(x)=\sqrt{\frac{x}{2\pi}},\\ 
& R(x)=\mcal{O}(x^{-1/4}),\\ 
& t\in [x,x+V] ,\ V\in[0,x^{1/4}], 
\end{split}
\edis  
where the functions 
\bdis 
\frac{2}{\sqrt{n}}\cos=\{t\omega_n(x)+\psi(x)\}
\edis 
are Riemann's oscillators with: 
\begin{itemize}
	\item[(a)] the amplitude 
	\bdis 
	\frac{2}{\sqrt{n}}, 
	\edis 
	\item[(b)] incoherent phase constant 
	\bdis 
	\psi(x)=-\frac x2-\frac{\pi}{8}, 
	\edis  
	\item[(c)] non-synchronized local times 
	\bdis 
	t\in [x,x+V], 
	\edis 
	\item[(d)] local spectrum of cyclic frequencies 
	\bdis 
	\{\omega_n(x)\}_{n\leq\tau(x)},\ \omega_n(x)=\ln\frac{\ln \tau(x)}{n}. 
	\edis 
\end{itemize}  

\begin{remark}
The Riemann-Siegel formula (see \cite{12}, p. 60) 
\bdis 
Z(t)=\sum_{n\leq \tau(x)}\frac{2}{\sqrt{n}}\cos=\{\vartheta(t)-t\ln n\}+\mcal{O}(t^{-1/4}) 
\edis 
is the Riemann's formula that has been restorted (by Riemann's manuscripts) and published by C.L. Siegel. 
\end{remark} 

\begin{remark}
Let us notice that by our opinion the Riemann-Siegel formula represents Riemann's fundamental contribution to the theory of oscillations (independently on the analytic number theory). Namely, the Riemann's oscillations are fated to describe of profound laws of our Universe. 
\end{remark} 

I would like to thank Micha Demetrian for his moral support of my study of Jacob's ladders.

\end{document}